\documentclass[a4paper,11pt,oneside,reqno]{amsart}
\usepackage[all,cmtip]{xy}
\usepackage{amsmath}
\usepackage{amsthm}
\usepackage{bbm} 
\usepackage{amssymb}
\usepackage{graphicx}
\theoremstyle{plain}
\newtheorem{teore}{Theorem}[section]
\newtheorem{defin}[teore]{Definition}
\newtheorem{lem}[teore]{Lemma}
\newtheorem{coro}[teore]{Corollary}
\newtheorem{propo}[teore]{Proposition}

\newtheorem{pro}{Problem}
\newtheorem{claim}{Claim}[teore]
\newtheorem*{claim*}{Claim}
\newtheorem*{thm*}{Theorem}
\newtheorem*{defi*}{Definition}
\theoremstyle{remark}
\newtheorem{ejemplo}[teore]{{\sc Example}}
\newtheorem{ejemplos}[teore]{{\sc Examples}}
\newtheorem{notas}[teore]{{\sc Remark}}
\newcommand{\nrm}[1]{\|#1\|}

\newcommand{\prop}{\begin{propo}}
\newcommand{\fprop}{\end{propo}}
\newcommand{\cor}{\begin{coro}}
\newcommand{\fcor}{\end{coro}}

\newcommand{\defi}{\begin{defin}\rm}
\newcommand{\fdefi}{\end{defin}}
\newcommand{\eje}{\begin{ejemplo}}
\newcommand{\feje}{\end{ejemplo}}
\newcommand{\ejes}{\begin{ejemplos}}
\newcommand{\fejes}{\end{ejemplos}}
\newcommand{\lema}{\begin{lem}}
\newcommand{\flema}{\end{lem}}
\newcommand{\teor}{\begin{teore}}
\newcommand{\fteor}{\end{teore}}
\newcommand{\nota}{\begin{notas}\rm}
\newcommand{\fnota}{ \end{notas}}
\newcommand{\clam}{\begin{claim}}
\newcommand{\fclam}{\end{claim}}
\newcommand{\clams}{\begin{claim*}}
\newcommand{\fclams}{\end{claim*}}

\newcommand{\lclam}{\begin{lclaim}}
\newcommand{\flclam}{\end{lclaim}}
\newcommand{\prucl}{\prue[Proof of Claim:]}
\newcommand{\fprucl}{\fprue}
\newcommand{\ben}{\begin{enumerate}}
\newcommand{\een}{\end{enumerate}}
\newcommand{\bit}{\begin{itemize}}
\newcommand{\eit}{\end{itemize}}

\newcommand{\mc}[1]{\mathcal{#1}}

\newcommand{\mk}[1]{\mathfrak{#1}}

\newcommand{\casos}{\begin{itemize}}
\newcommand{\fcasos}{\end{itemize}\setcounter{cs}{1}}

\newcommand{\fin}{\textsc{FIN}}

\newcommand{\ro}{\varrho}

\newcommand{\conj}[2]{ \{ {#1}\,:\,{#2} \} }

\newcommand{\ou}{\omega_{1}}

\newcommand{\om}{\omega}

\newcommand{\ip}{\sqsubseteq}

\newcommand{\ka}{\kappa}

\newcommand{\buit}{\emptyset}

\newcommand{\ga}{\gamma}
\newcommand{\id}{\text{Id }}

\newcommand{\Ga}{\Gamma}

\newcommand{\al}{\alpha}
\newcommand{\be}{\beta}
\newcommand{\de}{\delta}
\newcommand{\De}{\Delta}
\newcommand{\la}{\lambda}

\newcommand{\sig}{\sigma}

\newcommand{\vep}{\varepsilon}

\newcommand{\R}{{\mathbb R}}

\newcommand{\N}{{\mathbb N}}

\newcommand{\rest}{\upharpoonright}

\newcommand{\supp}{\mathrm{supp\, }}

\newcommand{\con}{\subseteq}

\newcommand{\cones}{\varsubsetneq }

\newcommand{\prue}{\begin{proof}}
\newcommand{\fprue}{\end{proof}}
\makeindex

\begin{document}
\title{Positional graphs and conditional structure of weakly null sequences}
\author{Jordi Lopez-Abad}
\author{Stevo Todorcevic}

\address{Instituto de Ciencias Matematicas (ICMAT).  CSIC-UAM-UC3M-UCM. C/ Nicol\'{a}s Cabrera 13-15, Campus Cantoblanco, UAM
28049 Madrid, Spain} \email{abad@icmat.es}

\address{Universit\'{e} Denis Diderot - Paris 7, C.N.R.S., UMR 7056, 2 place Jussieu
- Case 7012, 72521 Paris Cedex 05, France}
\address{and}
\address{Department of Mathematics, University of Toronto, Toronto, Canada, M5S 2E4 }
\email{stevo@math.toronto.edu} 


\begin{abstract}
We prove that, unless assuming additional set theoretical axioms, there are no reflexive space without unconditional sequences of density the continuum.  We give   for every integer $n$ there are
normalized weakly-null sequences of length $\om_n$ without unconditional subsequences. This together with a result of \cite{Do-Lo-To} shows that $\om_\omega$ is the minimal cardinal $\kappa$ that
could possibly have the property that every weakly null $\kappa$-sequence has an infinite unconditional basic subsequence . We also prove that for every cardinal number $\ka$ which is smaller than
the first $\om$-Erd\"{o}s cardinal there is a normalized weakly-null sequence without subsymmetric subsequences. Finally, we prove that  mixed Tsirelson spaces of uncountable densities must always
contain isomorphic copies of either $c_0$ or $\ell_p$, with $p\ge 1$.
\end{abstract}
\maketitle
\section{Introduction}

Maurey and Rosenthal \cite{Ma-Ro} were the first to construct a weakly null sequence $(x_n)_{n<\om}$ without
infinite unconditional subsequences. Their machinery of coding block-sequences of finite subsets of $\omega$
by integers and using this to define special functionals has become a standard tool in this area for
constructing examples of separable Banach spaces with conditional norms. The most famous of these is the
example of Gowers and Maurey \cite{Ga-Ma} of a separable reflexive space without unconditional basic
sequences. In \cite{Ar-Lo-To} a new coding of finite block-sequences of finite subsets of $\omega_1$ was
employed in building a reflexive space of density $\omega_1$ with no infinite unconditional basic
subsequence. While at that time it was not clear if the construction can be stretched up to $\omega_2,$
$\omega_3,$ etc, the paper \cite{Do-Lo-To} showed that there could be no such construction for density
$\om_\omega,$ at least if one is not willing to go beyond the standard axioms of set theory.  More precisely,
\cite{Do-Lo-To} shows that it is consistent that every weakly null sequence of length $\om_\omega$ must
contain an infinite unconditional basic subsequence. In the present paper we supplement this by extending the
Maurey-Rosenthal construction and  showing that for every integer $n$ there are normalized weakly-null
sequences of length $\om_n$ without infinite unconditional basic subsequences. This is achieved by a deeper
analysis of the combinatorial properties of families of finite subsets of $\om_n$ that could be used in
coding the conditionality inside a norm on $c_{00}(\omega_n).$ We express this using the notion of countable
chromaticity  for certain positional graphs on $[\om_n]^{<\om},$ a notion that seems of independent interest
and that might have other uses.

Recall that Argyros and Tolias \cite{Ar-To} have constructed another example of a non separable Banach space
without infinite unconditional basic sequences. Their space has density continuum and  is a dual of a
separable hereditarily indecomposable space and so highly non-reflexive. We explain this by showing that it
is consistent that every weakly null sequence of length continuum (and so, in particular, every reflexive
space of density continuum) must contain an infinite  unconditional basic subsequence. This is done by
connecting the problem with other classical combinatorial ideas such as, for example, free subsets of
algebraic structures, or  the possibility of extending the Lebesgue measure to all sets of reals, known as
the Banach problem.

Recall also that subsymmetric sequences $(x_n)_n$ (see Section \ref{notion}) are those that for any
$m_0<\dots<m_k$ and $n_0<\dots<n_k$ the linear extension of $x_{m_i}\mapsto x_{n_i}$ defines an isomorphism
of norm uniformly bounded.  It is a classical result that every non-trivial weakly-null  subsymmetric
sequence has an unconditional subsequence.  So, modulo an application of Rosenthal's $\ell_1$-theorem,
getting a subsymmetric basic sequence is just a step away from getting an unconditional one. While analyzing
the corresponding cardinal conditions on densities that guarantee the existence of each such sequence, we
have discovered the huge difference between these two notions. For example, we show that while a relatively
simple argument of Ketonen \cite{Ke} shows that every space of density at least an $\om$-Erd\"os cardinal
must contain an infinite sub-symmetric sequence, the converse is also true: For every cardinal $\kappa$ which
is not $\om$-Erd\"os cardinal there is a weakly null $\kappa$-sequence without infinite subsymmetric
subsequences.

As it is well known in this area,   the Tsirelson space \cite{Ts} was the first example of a space without
subsymmetric sequences. In the final section of our paper we explain this construction when lifted into the
realm of non separable spaces. For example, we show that so-called mixed Tsirelson spaces of uncountable
densities, unlike their separable counterparts, must always contain subsymmetric sequences, indeed
isomorphic copies of either $c_0$ or $\ell_p$, for some $p\ge 1$.

\section{Preliminaries and notation}\label{notion}

Given a set $X$, and a cardinal number $\la$, let $[X]^\la$ denote the collection of all   subsets of $X$ of
cardinality $\la$. Let $[X]^{\le \la}$ and $[X]^{<\la}$ denote the collection of subsets of $X$ of
cardinality $\le \la$ and $<\la$, respectively.

Let $\ka$ be a cardinal number. Given $s,t\in [\ka]^{<\om} $, we write $s<t$ to denote that $\max s<\min t$.
Given an integer $n$, let $\mathrm{Bs}_n(\ka)$ denote the collection of  all \emph{block sequences} of length
$n$ of finite subsets of $\ka$, i.e. all sequences $(s_i)_{i<n}$ of finite subsets $s_i\con \ka$ such that
$s_{i}<s_{i+1}$ for every $i<n-1$.  Let $\mathrm{Bs}_{<\om}(\ka)$ denote the set of all finite  block
sequences of finite subsets of $\ka$.

A family $\mc B$ of subsets of  $\ka$ is \emph{large} when $\mc B\cap [A]^n\neq\buit$ for every infinite
subset $A$ of $\ka$ and every $n<\om$. The family $\mc B$  is \emph{very-large }when for every infinite
subset $A$ of $\ka$ there is some infinite subset $B$ of $A$ such that $[B]^{<\om}\con \mc B$.

Recall that the order-type $\mathrm{otp}(W,R)$ of a well ordered set $(W,R)$ is the unique ordinal number
$\al$ such that there is an order-preserving bijection between $(W,R)$ and $\al$ when $\al$ is endowed with
its natural order. Given $A,B \con \ka $ of the same order-type, let $\vartheta_{A,B}:A\to B$ be the unique
order-preserving bijection between $A$ and $B$. Given $I\con \mathrm{otp}(A)$, let
$s[I]:=\vartheta_{\mathrm{otp}(A),A}"I$, and given $i\in \mathrm{otp}(A)$, let $(A)_i:=A[\{i\}]$.


Recall that a seminormalized sequence $(x_\al)_{\al<\ga}$ in a Banach space $(X,\nrm{\cdot})$  indexed in an
ordinal number $\ga$ is called a (Schauder) basic sequence when there is a constant $K\ge 1$ such that
$$\nrm{\sum_{\al<\la}a_\al x_\al}\le
K\nrm{\sum_{\al<\ga}a_\al x_\al}$$
 for every $\la<\ga$ and every sequence of scalars $(a_\al)_{\al<\ga}$.

A basic sequence $(x_\al)_{\al<\ga}$ is called (suppression) unconditional when there is a constant $K\ge 1$
such that
$$\nrm{\sum_{\al\in A}a_\al x_\al}\le
K\nrm{\sum_{\al<\ga}a_\al x_\al}$$
 for every $A\con\ga$ and every sequence of scalars $(a_\al)_{\al<\ga}$.

The sequence $(x_\al)_{\al<\ga}$ is called subsymmetric when there is a constant $K\ge 1$ such that for every
$s,t\con \ga$ of the same cardinality and every sequence $(a_n)_{n<|s|}$ of scalars have that
$$\frac{1}{K}\nrm{\sum_{n<|s|}a_n x_{\vartheta_{|s|,s}(n)}}\le \nrm{\sum_{n<|s|}a_n x_{\vartheta_{|t|,t}(n)}}\le K \nrm{\sum_{n<|s|}a_n x_{\vartheta_{|s|,s}(n)}}. $$

\section{Chromatic numbers of positional graphs}
In this section we define a positional graph whose vertex set is a
sufficiently large family of finite subset of some cardinal
$\kappa$ and show that a good coloring of this graph with
countably many colors can be used to define a normed space with a
weakly null $\kappa$-sequence with no infinite unconditional
subsequence.
\defi
Let $\mk G=(V,E)$ be a graph. The chromatic number $\chi(\mk G)$ is the minimal cardinal number $\ka$ such
that there is a coloring $c:V\to \ka$, called a \emph{good coloring of $\mk G$}, of the set of vertexes $V$
into $\ka$-many colors such that no two vertexes $v_0\neq v_1$ in an edge in $E$ have the same color $c$.
\fdefi
\eje \label{rjirjeirjexgfd} Let $\mk G_{\mathrm{card}}(\ka):=(V,E_\mathrm{card})$ be the graph with  vertexes
$V$ the set of all finite block sequences $(s_i)_{i<d}$, and let $E$ be the set of pairs
$((s_i)_{i<d},(t_i)_{i<l})$ such that
$$(|s_i|)_{i<d} \neq (|t_i|)_{i<l}.$$
Then $\mk G_{\mathrm{card}}(\ka)$ has countable chromatic number:  Let $d:\om^{<\om}\to\om$ be any 1-1
function, and for a finite block sequence $(s_i)_{i<d}$ in $\ka$ define $c((s_i)_{i<d}):=d((|s_i|)_{i<d})$.
Clearly $c$ is a good coloring.  \feje

Recall that for $A,B\con \ka$ we write $A\ip B$ when $A$ is initial part of $B$, i.e. when $A\con B$ and
$A<B\setminus A$.

\defi
Given an integer $n$ and two subsets $A$ and $B$ of $\ka$, we say that $A$ and $B$ are in $n$-$\De$-
\emph{position} when there are $I,J$ subsets of $A\cap B $ such that
\begin{enumerate}
\item[(a)] $I<J$, $A\cap B=I\cup J$ and $|J| \le n$.
\item[(b)] $I\ip A, B$.
\end{enumerate}
$A$ and $B$ are in \emph{$\De$-position} when they are in $0$-$\De$-position. Let $\mk G_n(\ka)$ be the graph
whose set of vertexes is $[\ka]^{<\om}$, and the set of edges $E_\mathrm{pos}$ is the family of pairs $(s,t)$
such that $s$ and $t$ are not in $n$-$\De$-position. Given a family $\mc B\con [\ka]^{<\om}$, we write $\mk
G_n(\mc B)$ to denote the restriction of the graph $\mk G_n(\ka)$ to set of vertexes $\mc B$.

\fdefi

It is easily seen (see, for example, \cite{To}, Chapter 3) that $\chi(\mk G_0(\omega))=\chi(\mk G_0(\ou))=\om$ but that
$\chi(\mk G_0(\omega_2))>\om.$ In fact, we have the following more general result

\teor\label{ijjijdffddf}
$\chi(\mk G_{n-2}(\omega_n))>\om$ for all $2\leq n<\omega$
\fteor

\begin{proof}
This follows from the standard fact (see, for example,  \cite{Di-To}) that for every
$c:[\omega_n]^n\rightarrow\omega$ there exist a sequence $s_0<s_1<\cdot\cdot\cdot<s_{n-1}$ of $2$-element
subsets of $\omega_n$ such that $c$ is constant on their product.
\end{proof}

The use of the polarized partition property in this proof suggest the following fact proved along the same lines.

\teor
If $\ka \rightarrow \left(\begin{array}{c}
2 \\
2\\
\vdots\\
\end{array}\right)^{<\om}_\lambda$  then $\chi(\mk G_{n}(\kappa))>\lambda$ for all $ n<\omega.$
\fteor

This partition property that has been analyzed in detail in \cite{Di-To} asserts that for every coloring
$c:[\kappa]^{<\omega}\rightarrow \lambda$ there is an infinite block sequence $(s_i)_{i<\om}$ of $2$-element
subsets of $\kappa$ such that $c$ is constant on $\prod_{i<n} s_i$ for all $n<\omega.$ Cardinals $\kappa$
with this partition property with $\lambda=2$ are called \emph{polarized cardinals}. It is know that (see,
for example, \cite{Di-To}) that every polarized cardinal satisfies the partition property for
$\lambda=2^{\aleph_0},$ so in particular, every polarized cardinal is greater than continuum. It is also
known (see \cite{Di-To}), that a polarized cardinal may not be greater than $2^{\aleph_1}$ and that, more
interestingly, that i is consistent that $\aleph_\omega$ is a polarized cardinal. Clearly if $\aleph_\omega$
is a polarized cardinal then it is the minimal polarized cardinals. The minimal polarized cardinals are
interesting because they satisfy the partition property for $\lambda=\theta^{\aleph_0}$ for every
$\theta<\kappa.$ So, in particular, the minimal polarized cardinal $\kappa$ has the property that
$\theta^{\aleph_0}<\kappa$ for all $\theta<\kappa$ (see \cite{Di-To}).  Their interest for us here comes from
the following immediate fact.

\cor
If $\kappa$ is the minimal polarized cardinal then $\chi(\mk G_{n}(\kappa))=\kappa$ for all $ n<\omega.$

\fcor

\cor It is consistent relative to the existence of a measurable
cardinal that for every $n<\om$ one has that $\chi(\mk G_{n}(\om_\omega))=\om_\omega$.
\fcor

The following  result shows that $\omega_\omega$ is in some sense
the minimal cardinal that could possibly have this property.

\teor\label{jiejrieide23fd}
For every $n<\om$ there is a very-large $\mc B_n\subseteq [\om_n]^{<\om}$ such that $\chi(\mk G_{2n-1}(\mc
B_n))=\om$.
\fteor
We postpone its proof to the Section \ref{nieedffdsfs}.

\section{Positional graphs and conditional norms}

The purpose of this section is to prove the following result that connects positional graphs with the existence of large weakly null sequences with no infinite unconditional basic subsequences.

\teor\label{lkojfjdfjdw2ere}
Suppose that for some cardinal $\ka$ and some $n<\om$ there is a very-large family $\mc
B\subseteq[\kappa]^{<\omega}$  such that $\chi(\mk G_n(\mc B))=\om.$  Then there is a norm $\|\cdot\|$ on
$c_{00}(\kappa)$ such that the sequence $(u_\gamma)_{\gamma<\kappa}$ of unit vectors of $c_{00}(\ka)$   is a
weakly null Schauder basis of the completion of $(c_{00}(\ka), \|\cdot\|)$ with no infinite unconditional
basic subsequences.

\fteor

\prue As with basically all known constructions of conditional
norms,  the basic trick of adding conditionality remains unchanged
since its first appearance  in the classical construction by
Maurey and Rosenthal \cite{Ma-Ro} of a weakly-null
$\omega$-sequence  without infinite unconditional subsequences.
 Let $n<\om$, let $\mc B$ be a very-large family on $\ka$ such that $\mk{G}_n(\mc B)=\om$ and let $c$ be a
good coloring of the graph  $(\mc B, E_\mathrm{card}\cup E_\mathrm{pos})$. By re-enumeration if needed, we
assume that $c$ takes values on a set $M\con \om\setminus n$ with the lacunary condition
\begin{equation}
\label{ir4jtiifjgfsfbf}
\sum_{m\in M}\sum_{l\in M\setminus\{m\}}\min\{\sqrt{\frac{l}{m}},\sqrt{\frac{m}{l}}\}\le 1,
\end{equation}
and such that $c(\buit)\ge n$. We say that a finite block sequence $(s_i)_{i<d}$ of subsets of $\ka$ is
\emph{$\mc B$-special} when
\begin{enumerate}
\item[(a)] $\bigcup_{j<i}s_j \in \mc B$  for every $i< d$.
\item[(b)] $|s_{i}|=c(\bigcup_{j< i}s_i )$ for every $i< d$.
\end{enumerate}
\clam
For every infinite set $A$ of $\ka$ and every $d$ there is a $\mc B$-special sequence of length $d$
consisting on subsets of $A$.
\fclam
\prucl
Let $B$ be an infinite subset of $A$ such that $[B]^{<\om}\con \mc B$. Now we simply choose any sequence
$(s_i)_{i<d}$  consisting of subsets of $B$ such that  (b) above holds.
\fprucl
Now let
$$K:=\conj{\sum_{i<d} \frac{1}{\sqrt{|s_i|}}\mathbbm{1}_{s_i}}{(s_i)_{i<d}\text{ is $\mc B$-special}}.$$
On $c_{00}(\ka)$ define the  norm $\nrm{\cdot}_K$ for $x\in c_{00}(\ka)$ by
$$\nrm{x}_K:=\max\{\nrm{x}_\infty,\sup\conj{\langle x,f   \rangle}{f\in K}\}.$$
Let $X$ be the completion of $(c_{00}(\ka),\nrm{\cdot}_K)$. It is not difficult to prove that
$(u_\ga)_{\ga<\ka}$ is a Schauder basis of $X$.
\clam
Suppose that $s\in [\ka]^{<\om}$ is such that $|s|\in M$. Then
\begin{equation}\label{mnfkkfmdfwere}
\nrm{\frac1{|s|^\frac{1}{2}}\mathbbm{1}_s}\le 2.
\end{equation}
Consequently,    $(u_\ga)_{\ga<\ka}$ is a weakly-null sequence.
\fclam
\prucl
Let $(s_i)_{i<d}$ be    $\mc B$-special with $d\le |s_0|$.  Then
\begin{align*}
\langle \frac{1}{|s|^\frac12}\mathbbm 1_s,\sum_{i<d} \frac{1}{|s_i|^\frac{1}2} \rangle\le  & \sum_{{i<d, |s_i|<|s|}}\frac{|s_i|^{\frac12}}{|s|^{\frac12}}+
\sum_{{i<d, |s_i|>|s|}}\frac{|s|^{\frac12}}{|s_i|^{\frac12}}+1 \le 2.
\end{align*}
Now if $(u_\ga)_\ga$ were not weakly null, we could find  $\vep>0$ and an infinite set $A\con \om_n$ such
that
\begin{equation}\label{ngjnnjdffgf}
\nrm{\frac{1}{|s|^\frac12}\mathbbm{1}_s}\ge \vep |s|^\frac12 \text{ for every $s\con A$ finite.}
\end{equation}
Clearly \eqref{ngjnnjdffgf} is in contradiction with \eqref{mnfkkfmdfwere}.
\fprucl
\clam
$(u_\ga)_{\ga<\ka}$ does not have infinite unconditional subsequences.
\fclam
\prucl
Fix a subset $A\con \ka$ of order type $\om$, and fix $L\ge 1$. We see that $(u_\ga)_{\ga\in A}$ is not
$L$-unconditional. Let $k<\om$ be such that $k>8L$. Let $x:=\sum_{i<k}(-1)^i\mathbbm{1}_{s_i}/{|s_i|^{1/2}}$,
where $(s_i)_{i<k}$ is a special sequence in $A$, and let
$y=\sum_{i<k/2}\mathbbm{1}_{s_{2i}}/{|s_{2i}|^{1/2}}$. Since $f:=\sum_{i<k}\mathbbm{1}_{s_i}/{|s_i|^{1/2}}\in
K$, it follows that $\nrm{y}_K\ge \langle f,y \rangle\ge k/2 $. We are going to see that $\nrm{x}_K\le 4$. To
this end,  fix $g=\sum_{i<d}\mathbbm{1}_{t_i}/{|t_i|^{1/2}}\in K$. Let
$$m_0:=\max\conj{i<\min\{d,k\}}{|s_i|=|t_i|}.$$
It follows then that $c(s_0,\dots,s_{m_0-1})=c(t_0,\dots,t_{m_0-1})$, and hence $|s_i|=|t_i|$ for all $i<m_0$
and, setting $s=\bigcup_{i<m_0}s_i$ and $t=\bigcup_{i<m_0}t_i$, there are $I,J\con s\cap t$ such that
\begin{enumerate}
\item[(a)] $I<J$, $s\cap t=I\cup J$ and $|J|< n$.
\item[(b)] $I\ip s, t$.
\end{enumerate}
Let $i_0<m_0$ be the last $i<m_0$ such that $s_i\con I$. Then it follows that
\begin{enumerate}
\item[(c)]$s_i=t_i$ for every $i\le i_0$, and
\item[(d)] $(\bigcup_{i=i_0+1}^{m_0-1}s_i)\cap (\bigcup_{i=i_0+1}^{m_0-1}t_i) $ has cardinality at most $n-1$.
\end{enumerate}
Since   $|s_i|=c(s_0,\dots,s_{i-1})$, $i<k$ and $|t_j|=c(t_0,\dots,t_{j-1})$, $j<d$, it follows that
$|s_i|\neq |t_j|$ for $i\neq j$. Hence,

\begin{align*}
|\langle g,x \rangle|\le  &
\left| \sum_{i<d}\sum_{j<k, |t_j|=|s_i|}\langle \frac{1}{\sqrt{|t_i|}}\mathbbm{1}_{t_i},\frac{(-1)^j}{\sqrt{|s_j|}}\mathbbm{1}_{s_j}\rangle\right|
+  \sum_{i<d}\sum_{j<k, |t_j|\neq |s_i|}\left|\langle \frac{1}{\sqrt{|t_i|}}\mathbbm{1}_{t_i},\frac{(-1)^j}{\sqrt{|s_j|}}\mathbbm{1}_{s_j}\rangle\right|
\le \\
\le & \left|\sum_{i\le i_0}\langle \frac{1}{\sqrt{|s_i|}}\mathbbm{1}_{s_i}, (-1)^i\frac{1}{\sqrt{|s_i|}}\mathbbm{1}_{s_i} \rangle\right|+
 \left| \sum_{i=i_0+1}^{m_0-1}\langle \frac{1}{\sqrt{|t_{i}|}}\mathbbm{1}_{t_{i}}, (-1)^i\frac{1}{\sqrt{|s_{i}|}}\mathbbm{1}_{s_{i}}  \rangle \right| + \\
 + & \left|\langle \frac{1}{\sqrt{|t_{m_0}|}}\mathbbm{1}_{t_{m_0}}, (-1)^i\frac{1}{\sqrt{|s_{m_0}|}}\mathbbm{1}_{s_{m_0}}   \rangle\right| +
  \sum_{i<d}\sum_{j<k, |t_j|\neq |s_i|}\left|\langle \frac{1}{\sqrt{|t_i|}}\mathbbm{1}_{t_i},\frac{(-1)^j}{\sqrt{|s_j|}}\mathbbm{1}_{s_j}\rangle\right| \le \\
  \le & 1+ \frac{\left|(\bigcup_{i=i_0+1}^{m_0-1}s_i)\cap (\bigcup_{i=i_0+1}^{m_0-1}t_i) \right|}{|s_{i_0}|}+ 2\le 4.
\end{align*}
\fprucl
\fprue
\nota
The previous example leads to a $c_0$-saturated space. It is possible to modify the construction to get a
reflexive example, tough unconditionally saturated.
\fnota

\begin{pro}
Is there a similar combinatorial condition on uncountable $\kappa$ that ensures the existence of a reflexive Banach space of density $\kappa$ with no infinite unconditional basic sequences?
\end{pro}

In \cite{Ar-Lo-To}, we have provided such an example for $\kappa=\omega_1.$

\section{Proof of Theorem \ref{jiejrieide23fd}}\label{nieedffdsfs}

Recall that the \emph{Shift graph}  on  a totally ordered set $(A, <)$ is the graph $([A]^{<\om}, \mathrm{Sf})$ whose   edges  are
the pairs $(s,t)$ of finite subsets of $\ka$ such that $s\setminus \{\min s\}\ip t$. The following  notions play an important role in the proof.

\defi
Let $A,X$ be two totally ordered sets, and let $f:[A]^n\to X$.
\begin{enumerate}
\item[(a)] We call $f$   \emph{Shift-increasing} if
\begin{equation}
\text{ $f(s)<f(t)$ for every $(s,t)\in \mathrm{Sf}$ in $A$.}
\end{equation}
\item[(b)] We call $f$ \emph{$\min$-dependant} if
\begin{equation}
\text{$f(s)=f(t)$ implies that $\min s=\min t$ for every $s,t\in [A]^n$.}
\end{equation}
\end{enumerate}
\fdefi

The proof crucially depends also on the following concept from \cite{To}.

\defi\label{rofunction}
A function $\ro:[\ka^+]^2\to \ka$  is called an (injective version of) $\ro$-function if
\begin{enumerate}
\item[(a)] $\ro$ is subbaditive, i.e. for every $\al<\be<\ga<\ka^+$
\begin{enumerate}
\item[(a.1)] $\ro(\al,\be)\le \max\{\ro(\al,\ga),\ro(\be,\ga)\}$,
\item[(a.2)] $\ro(\al,\ga)\le \max\{\ro(\al,\be),\ro(\al,\ga)\}$.
\end{enumerate}
\item[(b)] $\ro(\al,\be)\neq \ro(\bar \al,\be)$ for every $\al \neq \bar \al<\be$.
\item[(c)] $\ro(\al,\be)\neq \ro(\be,\ga)$ for every $\al<\be<\ga$.
\end{enumerate}
\fdefi

It is proved in \cite{To} (see Definition 3.2.1, Lemma 3.2.2 dealing with the case $\kappa=\omega$ and Chapter 9 for the general version) that such a function  $\ro:[\ka^+]^2\to \ka$ exists for every regular cardinal $\kappa.$

\defi
For each integer $n$ we fix an injective $\ro$ function $\ro^{(n)}$ on $\om_n$.
Let $n\in \om$. For each $i\le n$ we  define recursively $f_i^{(n)}:[\om_n]^{i+1}\to \om_{n-i}$ as follows:
\begin{enumerate}
\item[(1)] $f_0^{(n)}:=\id_{\om_n}$;
\item[(2)] $f_i(\al_0,\al_1,\dots,\al_{i}):=
\ro^{(n-(i-1))}(f_{i-1}(\al_0,\dots,\al_{i-1}),f_{i-1}(\al_1,\dots,\al_i))$ for each $\al_0<\dots<\al_i$
in $\om_n$ and each $0<i\le n$.
\end{enumerate}
Let $f_n:=f_n^{(n)}:[\om_n]^{n+1}\to \om$.
\fdefi
%
%

\prop\label{jbbeiu}
Suppose that $\al,\bar \al<\al_0<\dots<\al_{i-1}$ are such that
\begin{enumerate}
\item[(a)] $f_{i}^{(n)}(\al,\al_0,\dots,\al_{i-1})=f_{i}^{(n)}(\bar \al,\al_0,\dots,\al_{i-1})$.
\item[(b)] $f_{j}^{(n)}(\al,\al_0,\dots,\al_{j-1}),f_{j}^{(n)}(\bar \al,\al_0,\dots,\al_{j-1})<f_{j}^{(n)}(\al_0,\dots,\al_{j})$ for every $j<
i$.
\end{enumerate}
Then $\al=\bar \al$.
\fprop
\prue
This in done by induction on $i\ge 0$. The case $i=0$ is trivial. Suppose that $i>0$. Then
\begin{align*}
&\ro^{(n-(i-1))}(f_{i-1}^{(n)}(\al,\al_0,\dots,\al_{i-2}),f_{i-1}^{(n)}(\al_0,\dots,\al_{i-1}))=  f_{i}^{(n)}(\al,\al_0,\dots,\al_{i-1})=\\
&= f_{i}^{(n)}(\bar \al,\al_0,\dots,\al_{i-1})= \ro^{(n-(i-1))}(f_{i-1}^{(n)}(\bar \al,\al_0,\dots,\al_{i-2}),f_{i-1}^{(n)}(\al_0,\dots,\al_{i-1})).
\end{align*}
By the hypothesis (b), it follows that
$$f_{i-1}^{(n)}(\al,\al_0,\dots,\al_{i-2}),f_{i-1}^{(n)}(\bar
\al,\al_0,\dots,\al_{i-2}) <f_{i-1}^{(n)}(\al_0,\dots,\al_{i-1})).$$  So, by the property (b) of $\ro^{(n-(i-1))}$ in Definition
\ref{rofunction}, we get that $f_{i-1}^{(n)}(\al,\al_0,\dots,\al_{i-2})=f_{i-1}^{(n)}(\bar
\al,\al_0,\dots,\al_{i-2})$, and by inductive hypothesis we obtain that $\al=\bar \al$.
\fprue

\lema\label{ijijigh}
For  every $A\con [\om_n]^\om$ there is $B\in[A]^\om$ such that $f_i^{(n)}\rest [B]^{i+1}$ is
Shift-increasing for every $i\le n$.
\flema
\prue
Let $c:[A]^{n+2}\to {}^{n+1}{3}$ be the coloring defined for each $\al_0<\dots<\al_{n+1}$ in $A$ by
\begin{equation}
(c(\al_0,\dots,\al_{n+1}))_i=\left\{
\begin{array}{ll}
0 &\text{if $f_i^{(n)}(\al_0,\dots,\al_i)<f_i^{(n)}(\al_1,\dots,\al_{i+1})$}\\
1 & \text{if $f_i^{(n)}(\al_0,\dots,\al_i)=f_i^{(n)}(\al_1,\dots,\al_{i+1})$}\\
2 & \text{if $f_i^{(n)}(\al_0,\dots,\al_i)>f_i^{(n)}(\al_1,\dots,\al_{i+1})$}
\end{array}
\right.
\end{equation}
for each $i\le n$. Using the Ramsey Theorem we can find $B\in [A]^\om$ which is $c$-monochromatic, with
constant value $(\vep_i)_{i\le n}\in  {}^{n+1}{3}$.
\clam
$\vep_i=0$ for every $i\le n$.
\fclam
It follows easily from the claim that $B$ fulfills the requirement of the statement in the Lemma. So, it
rests to show the claim.
\prucl
We prove first that $\vep_i<2$ for  every $i\le n$. Otherwise, fix $i\le n$ such that $\vep_i=2$, and  let
$(s_n)_{n\in \om}$ be an arbitrary sequence in $[B]^{n+2}$ such that for every $n$ one has that
$(s_n,s_{n+1})\in \mathrm{Sf}_{n+2}$. For each $n$, let $t_n:=s_n[i+1]$. It follows then that
$f_i^{(n)}(t_n)>f_i^{(n)}(t_{n+1})$ for every $n$. This is impossible because the range of $f_i^{(n)}$ is the
well founded set $\om_{n-i}$.

Next, we prove that $\vep_i\neq 1$ for every $i\le n$, by induction on $i\le n$.  For $i=0$ we have that
$f_{0}^{(n)}(\al)=\al$ so it is clear that $\vep_0=0$.  Now suppose that $i>0$. Let $\al_0<\dots<\al_{i+1}$
in $B$. By inductive hypothesis,
\begin{equation}\label{io43hrih4t4}
f_{i-1}^{(n)}(\al_0,\dots,\al_{i-1})<f_{i-1}^{(n)}(\al_{1},\dots,\al_{i})<f_{i-1}^{(n)}(\al_{2},\dots,\al_{i+1}).
\end{equation}
It follows then by the property (c) of $\ro^{(n-(i-1))}$ that
\begin{align*}
f_{i}^{(n)}(\al_0,\dots,\al_{i})= & \ro^{(n-(i-1))}(f_{i-1}^{(n)}(\al_0,\dots,\al_{i-1}),f_{i-1}^{(n)}(\al_{1},\dots,\al_{i}))\neq \\
\neq &
 \ro^{(n-(i-1))}(f_{i-1}^{(n)}(\al_1,\dots,\al_{i}),f_{i-1}^{(n)}(\al_{2},\dots,\al_{i+1}))=f_{i}^{(n)}(\al_1,\dots,\al_{i+1}),
\end{align*}
so $\vep_i\neq 1$.
\fprucl
\fprue

\lema\label{mindep}
For  every $A\con [\om_n]^\om$ there is $B\in[A]^\om$ such that  the restriction $f_i^{(n)}\rest [B]^{i+1}$
is $\min$-dependant for every $i\le n$.
\flema
\prue
We use first Lemma \ref{ijijigh} to find $C\in [A]^\om$ such that $f_i^{(n)}\rest [C]^{i+1}$ are
Shift-increasing for every $i\le n$.  We use now the Erd\"{o}s-Rado Canonization  Theorem to find and infinite
set $B\in [C]^\om$  and for each $i\le n$ sets $J_i\con i+1 $ such that
\begin{equation}\label{dfksdfjkd2ww}
\text{$f_i^{(n)}(s)=f_i^{(n)}(t)$ if and only if $s[J_i]=t[J_i]$ for every $s,t\in [B]^{i+1}$ and every $i\le n$}.
\end{equation}
\clam
$0\in J_i$ for every $i\le n$.
\fclam
It is clear that this claim proves that $C$ has the desired properties.
\prucl
Fix $i\le n$. Since  \eqref{dfksdfjkd2ww} holds, it suffices to prove that if
$f_{i}^{(n)}(\al_0,\al_1,\dots,\al_{i})=f_{i}^{(n)}(\bar \al_0,\al_1,\dots,\al_{i})$, then $\al_0=\bar\al_0$.
Since $f_j^{(n)}\rest [B]^{j+1}$ are all Shift-increasing, we can apply Proposition \ref{jbbeiu} and get that
$\al_0=\bar \al_0$, so $0\in J_i$.
\fprucl
\fprue

\defi
Let $\mc B_n$ be the set of  all finite sets $s$ of $\om_n$ such that
\begin{enumerate}
\item[(a)] $f_n\rest [s]^{n+1}$ is $\min$-preserving.
\item[(b)] $f_i^{(n)}\rest  [s]^{i+1}$ is shift-increasing for every $i<n$.
\end{enumerate}
Let $c_n:\mc B_n\to \bigcup_{k<\om} {}^{k^{n+1}} \om $ be the coloring
$$c_n(s):=f_n\circ \vartheta_{|s|,s}^{n+1}.$$
\fdefi
\prop
$\mc B_n$ is a very-large family of finite subsets of $\omega_n$ and  $c_n$ is a good coloring of the graph $\mk G_n(\mc B_n)$.
\fprop
\prue
The Lemma \ref{ijijigh} and Lemma \ref{mindep} proves that $\mc B$ is very-large. Let us see that $c_n$ is a
good coloring.  Suppose that $s, t\in \mc B$ are such that $c_n(  s)=c_n(  t)$. If $|s\cap t|< 2n$, then $s$
and $t$ are in $n$-$\De$-position. So, suppose that $|s\cap t|\ge 2n$. Let
$J:=\{\ga_0<\dots<\ga_{n-1}<\ga_n<\dots<\ga_{2n-1}\}$ be the set of the last $2n$ elements of $s\cap t$. We
prove first that $\vartheta(\ga_i)=\ga_i$ for every $i<n$: Fix such $i<n$; then
\begin{equation}\label{jijjfjgfgfss}
f_n(\ga_i,\ga_n,\dots,\ga_{2n-1})=f_n(\vartheta_{s,t}(\ga_i),\vartheta_{s,t}(\ga_n),\dots,\vartheta_{s,t}(\ga_{2n-1}))
\end{equation}
Since $\vartheta_{s,t}"(s\cap t)\cup t$, and since $f_n\rest [t]^{n+1}$ is $\min$-dependant, it follows from
\eqref{jijjfjgfgfss} that $\vartheta_{s,t}(\ga_i)=\ga_i$.

\clam
$s\cap \ga_0 =t\cap \ga_0=s\cap t\cap \ga_0$.
\fclam
It is clear that the previous claim gives that $s$ and $t$ are in $n$-$\De$-position.
\prucl
Let $\ga\in s\cap \ga_0$. Then
\begin{equation}
f_n(\ga,\ga_0,\dots,\ga_{n-1})=f_n(\vartheta_{s,t}(\ga),\vartheta_{s,t}(\ga_0),\dots,\vartheta_{s,t}(\ga_{n-1}))=f_n(\vartheta_{s,t}(\ga),\ga_0,\dots,\ga_{n-1}).
\end{equation}
Since $f_i^{(n)}\rest [s]^{i+1}$ and $f_i^{(n)}\rest [t]^{i+1}$ are shift-increasing for every $i<n$, it
follows from Proposition \ref{jbbeiu} that $\vartheta_{s,t}(\ga)=\ga$, so $\ga\in t$. Similarly one proves
that $t\cap \ga_0\con s\cap \ga_0$.
\fprucl
\fprue

\section{Weakly-null sequences and the continuum}
Recall that at the present   there exist only two examples of non separable spaces with no infinite unconditional basic sequences, the example of   \cite{Ar-To} that has  of density $2^{\aleph_0}$
and the example of \cite{Ar-Lo-To} that has density $\aleph_1.$ While the space of \cite{Ar-To}  has density that it is at least consistently larger there is a crucial differences between these two
examples. The space of \cite{Ar-Lo-To} is reflexive while the space of \cite{Ar-To} is far from this, as it is a dual of a separable hereditarily indecomposable space. In this section we use some
known combinatorial properties of cardinals to show that this difference is indeed essential.
 More precisely, we show that there could be no constructions of weakly null sequences of length continuum with no unconditional
subsequences if one is not willing to use additional set theoretic assumptions such as, for example, the Continuum Hypothesis.
This follows from the following result.

\teor\label{mijtjijgff}
 It is consistent relative to the consistency of the existence of an $\omega$-Erd\"os cardinal\footnote{See the next section for definition.} that every weakly null sequence of length
 continuum contains an infinite unconditional basic sequence.
 \fteor
  Our proof will also reveals  the following interesting connection with the classical problem of Banach
  about extending the Lebesgue measure to all sets of reals.
 \teor\label{werdosdfweref}
Suppose that the Lebesgue measure extends to a total countably additive measure on $\R$.  Then every weakly null sequence of length
 continuum contains an infinite unconditional basic sequence.
 \fteor

 Our analysis is based on the following classical concept.
\defi
A cardinal $\kappa$ has the $\omega$-\emph{free-set property} if every algebra $\mc A$ on domain $A$ of cardinality $\kappa$ with no more than countably many operations (but with no restriction on
their arities) has an infinite \emph{free set}, an infinite subset $X$ of $A$ such that no $x\in X$ is in the subalgebra generated by $X\setminus\{x\}$. We use ${\rm Fr}(\kappa,\omega)$ to denote
this property of a cardinal $\kappa.$ \fdefi

It is known that every polarized cardinal $\kappa$  has the free-set property ${\rm Fr}(\kappa,\omega)$  (see \cite{Di-To}) and that it is this property that is more closely tied with the problem of
finding infinite unconditional subsequences of  a given  weakly null $\kappa$-sequence. More precisely, we have the following result from \cite{Do-Lo-To}.

\teor\label{free}
If  ${\rm Fr}(\kappa,\omega)$ holds then every weakly null $\kappa$-sequence contains an infinite unconditional basic subsequence.
\fteor

Unlike the first polarized cardinal, the first cardinal satisfying
${\rm Fr}(\kappa,\omega)$ does not need to be larger than the
continuum. There are several ways one can see this, the simplest
being the following fact which must have been already observed
before but we were not able to find a reference.

 \lema\label{ccc}
 The property ${\rm Fr}(\kappa,\omega)$ is preserved under forcing by the posets satisfying the countable chain condition.
 \flema

 \prue
 Given a poset $\mc P$ and a $\mc P$-name $\dot{\mc A}$ for an algebra on $\kappa,$ we may assume that the operations of $\dot{\mc A}$ are all coded up in a single function of the form $\dot{f}:\kappa^{<\om}\rightarrow\kappa.$ Since $\mc P$ satisfies the countable chain condition, for every name $\tau$ for an ordinal $<\kappa$ there is a countable set $C(\tau)\subseteq\kappa$ such that every $p\in\mc P$
 forces $\tau$ is an element of $C(\tau).$ So we can find a sequence $g_n:\kappa^{<\omega}\rightarrow\kappa$ $(n<\om)$ with the property that for
 every $s\in \kappa^{<\om}$
 $$C(\dot{f}(s))=\{g_n(s):n<\om\}.$$ Applying our assumption ${\rm Fr}(\kappa,\omega)$ to the algebra $(\kappa, g_n)_{n<\om}$ we get an infinite subset $X$ of $\kappa$ that is free
 in this algebra. It follows that every $p\in\mc P$ forces that $X$ is also free relative to the algebra $\dot{\mc A}$ in the forcing extension of $\mc P.$
 \fprue

 \prue[Proof of Theorem \ref{mijtjijgff}]
 Every $\omega$-Erd\"os cardinal has the property  ${\rm Fr}(\kappa,\omega)$ (see, for example, \cite{De}). Now the conclusion follows from Theorem \ref{free} and Lemma \ref{ccc}, by going to the
 forcing extension were at least $\ka$ many real numbers are added by a  poset satisfying the countable chain condition.
 \fprue

 \prue[Proof of Theorem \ref{werdosdfweref}]
 This follow by combining Theorem \ref{free} with the classical results of Erd\"os and Hajnal \cite{Er-Ha} and Solovay \cite{So}.
 More precisely, by the result of \cite{So}, the assumption implies the existence of a regular cardinal $\kappa\leq 2^{\aleph_0}$ with the
 \emph{Jonsson property}, i.e., with the property that every algebra with domain $\kappa$ contains a proper subalgebra of the same cardinality. On the other hand, it is proved in \cite{Er-Ha} (see also \cite{De}) that every Jonsson cardinal $\kappa$ has the free-set property   ${\rm Fr}(\kappa,\omega).$
 \fprue

 One of the reasons for mentioning Theorem \ref{werdosdfweref}
  is that it suggests using measure theoretic conditions on norms in order to have infinite unconditional subsequences of a given  large weakly null sequence.
   This, of course, remains to be explored, but we mention now one of such possible result.

\prop
Suppose that  $(x_\al)_{\al<\ka}$ be a normalized weakly null sequence \footnote{Indeed it suffices to assume
that the sequence is weakly countably-null} in some normed space $X.$ Suppose that the index-set $\kappa$
suppers a countably additive probability measure $\mu$ that gives measure zero to all countable subsets of
$\kappa$ and is defined on some $\sigma$-field of subsets of $\kappa$. Suppose that all norm-configurations
induced by subsets of the finite power $\kappa^n$  are measurable relative to the power measure $\mu^n.$ Then
$(x_\al)_{\al<\ka}$ contains an infinite unconditional basic subsequence.
\fprop

\prue
 For each finite set $s\con \ka$, choose a countable subset $N_s\con
S_{X^*}$ such that $\nrm{x}=\sup\conj{f(x)}{f\in N_s}$ for every $x\in \langle x_\al \rangle_{\al\in s}$. Let
$\theta:\fin(\ka)\to [\ka]^{\om}$ be defined for $s\in \fin(\ka)$ by
$$\theta(s):=\conj{\al<\ka}{\text{there is some $f\in N_s$ such that $f(x_\al)\neq 0$}}.$$
A finite set $s\in \fin(\ka)$ is called \emph{$\theta$-free} if
$$f(t)\cap s\con t\text{ for every $t\con s$}.$$
Let $F_n\con [\ka]^n$ be the set of $\theta$-free sequences of cardinality $n$. We identify a finite set
$s\in \fin(\ka)$ with the strictly increasing enumeration sequence $\vartheta_{|s|,s}\in \ka^n$.
\clam\label{kjjerjngjifdjg1}
Let $\mu$ be a $\sig$-additive  measure on $\ka$, and let $\mu_n$ denote the product measure of $\mu$ on
$\ka^n$. If the set $F_n$ is $\mu_n$-measurable, then $\mu_n([\ka]^n\setminus F_n)=0$.
\fclam
\prucl
Suppose otherwise that $\mu_n([\ka]^n\setminus F_n)>0$. For each $I\cones n$, and $j\in n\setminus I$, let
$$S_{I,j}:=\conj{s\in [\ka]^n}{(s)_j\in f(s[I])}.$$
It readily follows from the definition that
\begin{equation}
[\ka]^n\setminus F_n=\bigcup_{I\cones n,\, j\in n\setminus I} S_{I,j}.
\end{equation}
Let $I_0\cones n$ and $j_0\in n\setminus I_0$ be such that
$$\mu_n(S_{I_0,j_0})>0.$$
Set $J=I_0\cup \{j_0\}$, and let $\pi_{J}:[\ka]^n\to [\ka]^{J}$ be the canonical projection
$\pi_{J}(s)=s[J]$. It follows that
$$\mu_{m}(\pi_J(S_{J}))>0,$$
where $m=|J|$. By Fubini's Theorem, the set of $t\in [\ka]^{m-1}$ such that $\mu((\pi_J(S_J))_{t})>0$,
 has $\mu_{m-1}$-positive measure, where $(\pi_J(S_J))_t:=\conj{\al<\ka}{t\cup \{\al\}\in \pi_J(S_{J})}$.
 Observe that
\begin{align*}
(\pi_{J}(S_J))_t=&\conj{\al<\ka}{\text{there is some $s\in S_{J}$ such that $s[I_0]=t$ and $(s)_{j_0}=\al$}}\con f(t),
\end{align*}
so $(\pi_{J}(S_J))_t$ is countable, and hence $\mu(\pi_{J}(S_J))_t=0$, a contradiction.
\fprucl

\clam\label{kjjerjngjifdjg2}
Suppose that $s\in F_n$. Then $(x_\al)_{\al\in s}$ is a 1-unconditional basic sequence.
\fclam
\prucl
Let $(a_\al)_{\al\in s}$ be a sequence of scalars and fix $t\con s$. Let $f\in N_t$ be such that
\begin{equation}
f(\sum_{\al\in t}a_\al x_\al)=\nrm{\sum_{\al\in t}a_\al x_\al}.
\end{equation}
Since $s$ is $\theta$-free, it follows that $f(t)\cap s\con t$. This means that for every $\al\in s\setminus
t$ one has that $f(x_\al)=0$. So,
\begin{equation}
\nrm{\sum_{\al\in s}a_\al x_\al}\ge f(\sum_{\al\in s}a_\al x_\al)=f(\sum_{\al\in t}a_\al x_\al)=\nrm{\sum_{\al\in t}a_\al x_\al}.
\end{equation}
\fprucl

Recall that our assumption is that there is a non-trivial $\sig$-additive real-valued probability measure $\mu$
on $\ka$ such that for every $n\in \N$, $t\in F_{n-1}$ and every $i<n$, the sets
$$A_t^{i}:=\conj{\al<\ka}{\text{there is some $u\in F_n$ such that $u(i)=\al$ and $u[n\setminus \{i\}]=t$}}$$
are all $\mu$-measurable. So what we have accomplished is that  the sets $F_n$ are $\mu_n$-measurable   for
every $n\in \N$. It follows from the Claim \ref{kjjerjngjifdjg1} and the Claim \ref{kjjerjngjifdjg2} that the
sets $U_n:=\conj{s\in [\ka]^n}{(x_\al)_{\al\in s}\text{ is $1$-unconditional}}$ have $\mu_n$-measure 1 for
every $n$. Let $\mu_\om$ denote the infinite product measure on $\ka^{\om}$. Then it follows that the set of
sequences $(\al_n)_{n\in \om}$ such that $(x_{\al_n})_{n\in \om}$ is 1-unconditional has measure 1, and  we
are done.
\fprue

\section{Subsymmetric sequences}

\defi
Let $A$ be a set of ordinal numbers and let $\mc B\con [A]^{<\om}$. The family $\mc B$ is pre-compact when
its topological closure consists only on finite sets, when we consider the topology product topology on
$2^A$. The family $\mc B$ is compact when $\mc B$ is closed with respect to the product topology on $2^A$.
The family $\mc B$ is hereditary when $t\con s\in \mc B$ implies that $t\in \mc B$.  Recall that we say that
the family $\mc B$ is \emph{large} in $A$ when for every $B\in [A]^\om$ and every $n\in \N$ we have that $\mc
B\cap [B]^n\neq \buit$.
\fdefi
It is easy to see that if $\mc B$ is hereditary, then $\mc B$ is pre-compact if and only if $\mc B$ is compact.  It is also an exercise to prove that $\mc B$ is pre-compact if and only if   for every
infinite sequence $(s_n)_{n\in \om}$ in $\mc B$ has a subsequence $(s_n)_{n\in B}$ forming a $\De$-system, i.e. there is some $s\in [A]^{<\om}$ such that for every $n\neq m$ in $A$ one has that
$s_n\cap s_m=s$.

\defi
Let $\ka$ be a cardinal number and  $c:[\ka]^{<\om}\to X$. A
subset $H$ of $\ka$ is called \emph{$c$-homogeneous} if for every
$n<\om$ the restriction $c\rest [H]^{n}$ is constant. Given such
$c$, we define $\mc H(c)$ as the family of all $c$-homogeneous
subsets of $\ka$.

Recall that a cardinal number $\ka$ is called $\om$-Erd\"{o}s if
$$\ka\longrightarrow (\om)^{<\om}_2,$$
i.e., for every $c:[\ka]^{<\om}\to 2$  the family $\mc H(c)$ consists only on finite subsets of $\ka$. \fdefi In general, $\mc H(c)$ is a compact and hereditary family of $\mc P(\ka)$. It is not
difficult to see that $\ka$ is $\om$-Erd\"{o}s if and only if  every $c:[\ka]^{<\om}\to 2^{\aleph_0}$  there are no infinite $c$-homogeneous subsets of $\ka$. \prop\label{nm34jnmt4tgrgvc} Let
$c:[\ka]^{<\om}\to 2$. The family $\mc H(c)$ is a large compact and hereditary family  of (not necessary finite) subsets of $\ka$.
\fprop

\prue It is clear that $\mc H(c)$ is hereditary and closed in $\mc
P(\ka)$. Let us prove that $\mc H(c)$ is large. Given integers
$k,l,m$, let $R(k,l,m)$ be the minimal integer such that whenever
$d:[R(k,l,m)]^{k}\to m$ there is a subset $A\con R(k,l,m)$ of
cardinality $l$ such that $d\rest [A]^k$ is constant.

Fix now an infinite subset $A\con \ka$ and $n<\om$. Let  $s$ be a finite subset of $A$ of cardinality $|s|\ge
R(n,2n-1,2^{n})$. Let  $d:[s]^{n}\to 2^{n}$ be defined for $t\in [s]^n$ by $d(t):= (c(t[i]))_{i<n}$. Let
$u\in [s]^{2n-1}$ be such that $c\rest [u]^{n}$ is constant, with value $(\vep_i)_{i<n}\in 2^{n}$. Let
$v:=u[n]$. We claim that $v$ is $c$-homogeneous:  This is an easy consequence of the fact that for every
$i<n$ and every $w\in [v]^i$ there is $z\in [u]^n$ such that $z[i]=w$.
\fprue

Note that $\om$ is not $\om$-Erd\"{o}s: Let  $\mc S$ denote the Schreier family on $\om$. Then the characteristic function $\mathbbm{1}_{\mc
 S}$ does not have infinite homogeneous sets:
Suppose that there is some infinite $A\con \om$ such that $\mathbbm{1}_\mc S\rest [A]^{n}$ is constant, with
value $\vep_n\in 2$ for every $n<\om$. Since $\mc S$ is hereditary, we get that
\begin{equation}
\vep_n=\min_{m<n}\vep_m \text{ for every $n<\om$.}
\end{equation}
This means that, either there is some $m$ such that $\vep_m=0$, and hence $\mc S\cap [A]^m=\buit$ contradicting the fact that $\mc S$ is large, or else $\vep_m=1$ for all $m$, hence $[A]^{<\om}\con
\mc S$, contradicting the compactness of $\mc S$. Following this idea, the general situation is completely understood.

\prop
Let $\ka$ be an infinite cardinal. The following are equivalent:
\begin{enumerate}
\item[(1)] $\ka$ is $\om$-Erd\"{o}s.
\item[(2)] Every separated normalized sequence $(x_\al)_{\al<\ka}$ has a subsymmetric subsequence.
\item[(3)] There are no large compact and hereditary families on $\ka$.
\end{enumerate}
\fprop
\prue
(1) implies (2): This is due to \cite{Ke} but for the convenience of the reader, we sketch the proof here. Suppose that $\ka$ is a $\om$-Erd\"{o}s
cardinal, and suppose that   $(x_\al)_{\al<\ka}$ is a normalized separated sequence.  Let $c:[\ka]^{<\om}\to
2^{\aleph_0}$ be defined as follows. Let $s\in [\ka]^{<\om}$.  Note that $\vartheta_{s,|s|}$ linearly extends
to $\vartheta_{s,|s|}:\langle x_\ga \rangle_{\ga\in s}\to \langle u_n \rangle_{n\in |s|}$ by
$\vartheta_{s,|s|}(u_\ga):=u_{\vartheta_{s,|s|}(\ga)}$ for each $\ga\in s$.  Let $c(s):= (|s|,  (f\circ
\vartheta_{|s|,s})_{f\in (\langle x_\ga \rangle_{\ga\in s})^*})$. Let $A\con \ka$ be an infinite set such
that $c\rest [A]^n$ is constant for every $n$. This implies that $(x_\ga)_{\ga\in A}$ is 1-subsymmetric.

(2) implies (3): Suppose that $\mc B$ is a large compact and hereditary family  on $\ka$. Now, we use the
family $\mc B$ to define the following Schreier norm in $c_{00}(\ka)$: For $x\in c_{00}(\ka)$ let
$$\nrm{x}_\mc B:=\max\conj{|\langle x,\chi_s \rangle|}{s\in \mc B}.$$
Observe that $[\ka]^1\con \mc B$, so the formula above defines a norm on $c_{00}(\ka)$.  Let $X_\mc B$ be the
completion of $c_{00}(\ka)$ with respect to the norm $\nrm{\cdot}_\mc B$. It readily follows from the fact that
$\mc B$ is hereditary,  that the unit Hamel basis $(u_\ga)_{\ga< \ka}$ is an unconditional Schauder basis of
$X_{\mc B}$.
\clam
$(u_\ga)_{\ga< \ka}$ is a normalized weakly-null sequence without subsymmetric subsequences.
\fclam
\prucl
 Fix $A \in [\ka]^{\om}$. We first check that $(u_\ga)_{\ga\in A}$ is weakly-null. This is a standard fact:
Suppose that there is some  $f\in (X_\mc B)^*$ and   $\vep>0$ the set $B:=\conj{\ga\in A}{f(u_\ga)\ge \vep}$ is infinite. Now since $\mc B\rest B$ is a compact family on $B$, there is by \emph{Ptak's
Lemma},   some $\mu\in S_{\ell_1(B)}^+\cap c_{00}(B)$ such that
\begin{equation}\label{jtjgijthyt}
\text{$|\langle \mu,\mathbbm{1}_s \rangle|\le \frac\vep 2$ for every $s\in \mc B\rest A$.}
\end{equation}
Then,
\begin{equation}
\frac\vep 2\ge \nrm{ \sum_{\ga\in \supp \mu} (\mu)_\ga u_\ga}_\mc B \ge f(\sum_{\ga\in \supp \mu} (\mu)_\ga u_\ga)\ge \vep,
\end{equation}
which is, of course, impossible.

Now we prove that $(u_{\ga})_{\ga\in A}$ is not subsymmetric: Fix $C\ge 1$.  Use again Ptak's Lemma to find $\mu\in S_{\ell_1(A)}^+\cap c_{00}(A)$ such that $|\langle \mu,\chi_s \rangle|\le C/2$ for
every $s\in \mc B\rest A$. Let as above  $x:=\sum_{\ga\in \supp \mu} (\mu)_\ga u_\ga$. Then $\nrm{x}_\mc B\le C/2$. Since $\mc B$ is large, there is some $s\in \mc B \rest [A]^{n}$, where $n=|\supp
\mu|$. Let $y:=\sum_{\ga \in s}(\mu)_{\vartheta_{s,\supp \mu}(\ga)} u_{\ga}$. Then
\begin{equation}
\nrm{y}_\mc B\ge \langle y, \chi_s \rangle=\sum_{\ga\in \supp \mu} (\mu)_\ga=1 > \frac{1}{C}\nrm{x}_\mc B.
\end{equation}
\fprucl

(3) implies (1): Suppose that $\ka$ is not $\om$-Erd\"{o}s. Fix a
coloring  $c:[\ka]^{<\om}\to 2$ without infinite $c$-homogeneous
sets. Then $\mc H(c)$ is, by Proposition \ref{nm34jnmt4tgrgvc}, a
large compact and hereditary family of finite subsets of $\ka$.
\fprue

%

\section{Non-separable spaces of Tsirelson type }

We have seen that for every non $\om$-Erd\"{o}s cardinal $\ka$ there is a large family $\mc B$ on $\ka$. For $\ka=\om$, these families can be used to provide a separable Banach space without subsymmetric
sequences. A large family on $\om$ was used   to build the Tsirelson space\footnote{indeed, its dual} which is a reflexive space with an unconditional basis and without subsymmetric basic sequences.
Its construction is naturally generalized to an arbitrary infinite cardinal number $\ka$. The goal of this subsection is to prove that non-separable Tsirelson spaces have always subsymmetric basic
sequences, indeed isomorphic copies of some of the classical sequence spaces $c_0$ or $\ell_1$, $p\ge 1$.

Towards this goal, in all this subsection we fix an infinite cardinal  number $\ka$, an  \emph{hereditary}
 family $\mc B$ of finite subsets of $\ka$, and a real number $0<\theta<1$.
\defi
We say that a finite block sequence  $(E_i)_{i<n}$ of finite subsets of $\ka$ is \emph{$\mc B$-admissible}
when there is $\{\ga_i\}_{i<n}\in \mc B$ such that   $\ga_0\le E_0<\ga_1\le E_1<\dots  <\ga_{n-1}\le E_n$.
Similarly, a finite sequence $(x_i)_{i=1}^n$ of vectors of $c_{00}(\ka)$ is called $\mc B$-admissible when
$(\supp x_i)_{i=1}^n$ is $\mc B$-admissible.

Given such $\mc B$ and a real number $0<\theta<1$, we define the \emph{Tsirelson-like space} $T_{\theta,\mc
B}:= (T_{\theta,\mc B},\nrm{\cdot}_{\theta,\mc B})$ on $\ka$ as follows: The norm $\nrm{\cdot}_{\theta,\mc B}$ is the
unique norm on $c_{00}(\ka)$ satisfying that for every $x\in c_{00}(\ka)$ we have that
\begin{equation}
\label{jkgfbjfgnf}
\nrm{x}_{\theta,\mc B}=\max\{\nrm{x}_\infty,\sup\conj{\theta\cdot\sum_{i<n} \nrm{E_i x}_{\theta,\mc B}}{ (E_i)_{i<n}\text{ is $\mc B$-admissible}}\}.
\end{equation}
Then $T_{\theta,\mc B}$ is the completion
of the normed space $(c_{00}(\ka),\nrm{\cdot}_{\theta,\mc B})$.
\fdefi
Since the family $\mc B$ is hereditary, it follows easily that
\begin{equation}
\label{nnfdngjf}\text{if $(E_i)_{i<n}$ is $\mc B$-admissible, and
$F_i\con E_i$, $i<n$, then so is $(F_i)_{i<n}$.}
\end{equation}
This fact readily implies that the unit basis $(u_\ga)_{\ga< \ka}$ is a 1-unconditional Schauder basis of the
Tsirelson-like space $T_{\theta,\mc B}$.

The classical Tsirelson example is $T:=T_{1/2,\mc S}$, where $\mc S$ is the Schreier family on $\om$. An interesting result of Bellenot \cite{Be} states that  $T_{\theta,[\ka]^{\le n}}$  is
isomorphic to $c_{0}(\ka)$, if $\theta n \le 1$ or to $\ell_p(\ka)$ for $p=\log(n)/(\log(n)+\log(\theta))$.  In general, the situation is well understood for a compact family $\mc B$ in $\om$: If the
Cantor-Bendixson rank of $\mc B $ is an integer $n$, then the space $T_{\theta,\mc B}$ is saturated by copies of either $c_0$ when $\theta n\le 1$, or by copies of $\ell_p$, with
$p=\log(n)/(\log(n)+\log(\theta))>1$ if $\theta n> 1$. In this latter case, the space $T_{\theta,\mc B}$ is therefore, by a classical result of James, reflexive (see \cite{Be-De}). If otherwise the
Cantor-Bendixson rank of $\mc B$ is infinite, then the space $T_{\theta,\mc B}$  does not contain subsymmetric sequences, and therefore reflexive (see for example \cite{Lo-Ma}). In particular,
$T_{\theta,\mc B}$ never contains a copy of $\ell_1$.   If $\ka>\om$, the situation is very much different. To visualize this difference easily, let us assume that $\ka$ is not a cardinal but the
ordinal number $\om^2$, and let $\mc B$ be the family
$$\mc B:=\{\buit\}\cup\conj{\{\om \cdot i +n\}_{i<m}}{n<\om, \text{ and $m\le n$}}.$$
Then the only accumulation point in $\mc B$ is $\buit$, and hence its Cantor-Bendixson rank is $1$. However
every sequence  $(u_{\ga})_{\ga\in s}$ with $s\con \{\om\cdot n\}_{n>0}$ is $\mc B$-admissible, and hence
$(u_{\om\cdot n})_{n>0}$ is $\theta^{-1}$-equivalent to the unit basis of $\ell_1$. It is easy to modify $\mc
B$ to make it with infinite Cantor-Bendixson rank and still having a copy of $\ell_1$.

The main result is the following.
\teor\label{io43oi4rgrg}\mbox{ }
\begin{enumerate}
\item[(1)] Suppose that $\ka$ is an uncountable regular cardinal. Then there is a club $C$ of $\ka$ such that
either
\begin{enumerate}
\item[(1.1)] the sequence $(u_\ga)_{\ga\in C}$ in $T_{\theta,\mc B}$ is $\theta^{-1}$-equivalent to the unit basis of $\ell_1(\ka)$
or
\item[(1.2)] the closed linear span of ${\langle u_\ga \rangle_{\ga\in C}}$ in $T_{\theta,\mc B}$ is either
$c_0$-saturated or $\ell_p$-saturated for some $p\ge 1$.
\end{enumerate}
\item[(2)] Suppose that $\ka$ is an uncountable singular cardinal. Then for every cardinal number $\la<\ka$ there is a subset $C_\la$ of $\ka$ of
cardinality $\la$ such that either
\begin{enumerate}
\item[(1.1)] the sequence $(u_\ga)_{\ga\in C_\la}$ in $T_{\theta,\mc B}$ is $\theta^{-1}$-equivalent to the unit basis of $\ell_1(\la)$
or
\item[(1.2)] the closed linear span of ${\langle u_\ga \rangle_{\ga\in C}}$ in $T_{\theta,\mc B}$ is either
$c_0$-saturated or $\ell_p$-saturated for some $p\ge 1$.
\end{enumerate}
\end{enumerate}
\fteor

\defi Let $n\in \N$.  The  \emph{$(\mc B,n)$-Namba game}, or simply the $(\mc B,n)$-game $\Game_{\mc B,n}
$ is the following game of height $n$:  The first Player $I$ starts with an ordinal $\ga_0<\ka$, and the
Player $II$ replies with an ordinal $\ga_0< \eta_0<\ka$. Then the Player $I$ replies with $\eta_0<\ga_1<\ka$,
and then Player $II$ chooses $\ga_1< \eta_1<\ka$, and so on. The game ends after $n$ runs of each player.
Player $I$ wins if the set $\{\ga_i\}_{i<n}\in \mc B$; otherwise $II$ wins.
\fdefi

\defi A \emph{strategy} for Player $I$ or Player $II$ in the game $\Game_{\mc B,n}$ is a mapping $\sig:[\ka]^{<n-1}\to \ka$
 such that $s<\sig(s)$ for every $s\in
[\ka]^{<n-1}$.  An strategy $\sig$ for Player $I$ is \emph{winning} for the game $\Game_{\mc B,n}$ if
$\{\sig(\buit),\sig(\{\eta_0\}),\sig(\{\eta_0,\eta_1\}),\dots,\sig(\{\eta_0,\eta_1,\dots,\eta_{n-2}\})\}$ is
in $\mc B$ for every $\eta_0<\eta_1<\dots<\eta_{n-2}<\ka$.  A strategy $\sig$ for Player $II$ is
\emph{winning} for the game $\Game_{\mc B,n}$ if $\{\ga_0,\ga_1,\dots,\ga_{n-1}\}$ is not in $\mc B$ for
every $\ga_0<\ga_1<\dots<\ga_{n-1}<\ka$ such that
$\ga_0<\sig(\{\ga_0\})<\ga_1<\sig(\{\ga_0,\ga_1\})<\dots<\ga_{n-2}<\sig(\{\ga_0,\dots,\ga_{n-2}\})<\ga_{n-1}$.
\fdefi
Few useful facts.
\prop
\begin{enumerate}
\item[(a)] The game $\Game_{\mc B,n}$ is determined, i.e., either Player $I$ has a winning strategy or Player
$II$ has a winning strategy.
\item[(b)] If player $I$ has a winning strategy for the game $\Game_{\mc B,n}$ then he also has a winning
strategy for the game $\Game_{\mc B,m}$ for every $m\le n$. Symmetrically, if player $II$ has a winning
strategy for the game $\Game_{\mc B,n}$ then he also has a winning strategy for the game $\Game_{\mc
B,m}$ for every $m\ge n$.
\end{enumerate}
\fprop
\prue
(a) is a consequence of the fact that the game $\Game_{\mc B,n}$ is finite, (b) and (c) are consequences of
the fact that $\mc B$ is hereditary.
\fprue
\defi
A set $C\con \ka$ is closed under a strategy $\sig$ of Player $I$ or Player $II$ for the game $\Game_{\mc
B,n}$ when $\sig(s)< C/\ga$ for every  $\ga\in C$ and  $s\in [\ga+1 ]^{<n-1}$.
\fdefi
\prop\label{klkmniuuirr}
Suppose that $\ka$ is an uncountable regular cardinal. Then for every strategy $\sig$ of Player $I$ or Player
$II$ there is a club $C\con \ka$ closed under $\sig$.
\fprop
\prue
Fix $\sig:[\ka]^{<n-1}\to \ka$. We define $C:=\{\ga_\xi\}_{\xi<\ka}$ inductively: If $\xi$ is a limit
ordinal, then $\ga_\xi:=\sup_{\eta<\xi} \ga_\eta$. For the successor case $\xi+1$, Let
$$\ga_{\xi+1}:=\sup \conj{\sig(s)}{s\in [\ga_\xi+1]^{<n-1}}+\om.$$
Since $\ka$ is uncountable and regular $\ga_\xi<\ga_{\xi+1}<\ka$. It is easy to see that  $C$ has the
required properties.
\fprue
\defi
Given a family $\mc B$ of finite subsets of $\ka$, let
$$\al(\mc B):=\sup \conj{n<\om}{\text{ Player $I$ has a winning strategy for the game $\Game_{\mc B, n}$}}.$$
\fdefi

\defi
Let $C\con \ka$.
\begin{enumerate}
\item[(a)] $C$ is \emph{$(\mc B,I)$-closed }when for every \emph{integer} $n\le \al(\om)$ $C$ is closed under a
winning strategy for the game $\Game_{\mc B,n}$.
\item[(b)] $C$ is \emph{$(\mc B,II)$-closed} when $C$ is closed under a winning strategy of Player $II$ in the game $\Game_{\mc B,\al(\mc
B)+1}$, if $\al(\mc B)$ is an integer.
\item[(c)] $C$ is \emph{$\mc B$-strategically closed} when it is $(\mc B,I)$-closed and $(\mc B,II)$-closed.
\end{enumerate}
\fdefi
\prop
\begin{enumerate}
\item[(a)] The three notions above are hereditary.
\item[(b)] Suppose that $\ka$ is an uncountable regular cardinal. Then there is a club $C$ of $\ka$ which is $\mc
B$-strategically closed.
\end{enumerate}
\fprop
\prue
The first part is trivial. The second one  is an immediate consequence of Proposition \ref{klkmniuuirr} and
the fact that a countable intersection of clubs of $\ka$ is also a club of $\ka$.
\fprue
Recall that our main goal is to find a large subset $C$ of $\ka$ and identify the closed subspace of
$T_{\theta,\mc B}$ spanned by $\{u_\ga\}_{\ga\in C}$.  This would be relatively easy to do if the sequence
$(u_\ga)_{\ga\in C}$ in $T_{\theta,\mc B}$ were equivalent to itself as a sequence in the space
$T_{\theta,\mc B\rest C}$, because we would have reduced the main work to study restrictions of the family
$\mc B$.  However this is not the case, but there is another family $\mc F$ of finite subsets of $C$ for
which the sequence $(u_\ga)_{\ga\in C}$ in $T_{\theta,\mc B}$ and in $T_{\theta,\mc F}$ are 1-equivalent (see
Proposition \ref{kloihreiumnnds}).
\defi
Given $\Ga\con \ka$, $\Ga\neq \buit$ and without maximum, let $\varpi_C:\sup C\to C$ be defined by
$\varpi(\ga)=\min \conj{\de\in C}{\de\ge \ga}$. It is clear that $\varpi_C$ is onto. Given a family $\mc B$
of finite subsets of $\ka$, let
$$\varpi_\Ga(\mc B):=\conj{\varpi_{\Ga}"(s)}{s\in\mc B}.$$
Let also
\begin{align*}
\Ga^+:=&\conj{\ga\in \Ga}{\text{there is some $\eta\in \Ga$ such that $]\de,\ga]=\{\ga\}$}}.
\end{align*}
In other words,  $\Ga^+$ is the set of successors of $\Ga$.
\fdefi
Observe that for $\Ga_0\con \Ga_1$, the mappings $\varpi_{\Ga_0}$ and $\varpi_{\Ga_1}\rest \sup \Ga_0$ are in
general different. However, it is easy to see that if $\Ga_0$ is an interval of $\Ga_1$, then the
corresponding two mappings coincide. It is also clear from the definition that the family $\varpi_{\Ga}(\mc
B)$ is hereditary because so is $\mc B$, and that $\mc B\rest C\con \varpi_{C}(\mc B)$.

The following two are useful facts.
\prop \label{mnkdwewrdfnknke}
Suppose that $(s_i)_{i<d}$ is a block sequence of subsets of $\Ga$ with $\{\min s_i\}_{i<d}\in
\varpi_{\Ga}(\mc B)$. Then the  sequence $(s_i)_{i<d}$ is $\mc B$-admissible.
\fprop
\prue
Let $\{\de_0<\dots<\de_{d-1}\}\in \mc B$ be such that $\min s_i=\varpi_{\Ga}(\de_i)$ for $i<d$.  Then, by
definition of $\varpi_{\Ga}$ and the fact that $\varpi_{\Ga}(\de_{i+1})=\ga_{i+1}>\max s_i$, we have that
$\ga_i<\de_{i+1}$ for every $i<d-1$. Hence, $\de_0\le \ga_0\le \max s_0<\de_1\le \ga_1\le \max
s_1<\dots<\de_{n-1}\le \ga_{n-1}$ and consequently $(s_i)_{i<d} $ is $\mc B$-admissible.
\fprue
\prop\label{kloihreiumnnds}
Let $x\in c_{00}(\Ga)$. Then $\nrm{x}_{\theta,\mc B}=\nrm{x}_{\theta,\varpi_{\Ga}(\mc B)}$.
\fprop
\prue
The proof is by induction on the cardinality $k$ of the support of $x$.  If $k=1$, the result is trivial. The
case $k>1$ readily follows from the following.
\clam
Let $(E_i)_{i<d}$ be a block sequence of finite subsets of $\Ga$. Then $(E_i)_{i<d}$ is $\mc B$-admissible if and only if $(E_i)_{i<d}$ is $\varpi_\Ga(\mc B)$-admissible.
\fclam
\prucl
Suppose that $(E_i)_{i<d}$ is  $\mc B$-admissible.   Let $\{\eta_i\}_{i<d}\in \mc B$ be such that $\eta_0\le
E_0<\eta_1\le E_1<\dots<\eta_{d-1}\le E_{d-1}$. Then $\varpi_{\Ga}(\eta_0)\le E_0<\varpi_{ \Ga}(\eta_1)\le
E_1<\dots<\varpi_{\Ga}(\eta_{d-1})<E_{d-1}$. Since $\{\varpi_{\Ga}(\eta_i)\}_{i<d}\in \varpi_{\Ga}(\mc B)$,
it follows that $(F_i)_{i<d}$ is also $\mc F$-admissible.

Suppose now that  $(E_i)_{i<d}$ is  $\varpi_{\Ga}(\mc B)$-admissible. Let $\{\ga_i\}_{i<d}\in
\varpi_{\Ga}(\mc B)$ be such that $\ga_0\le E_0<\ga_1\le E_1<\dots<\ga_{d-1}\le E_{d-1}$, and let
$\{\eta_i\}_{i<d}\in \mc B$ be such that $\ga_i:=\varpi_{\Ga}(\eta_i)$, $i<d$.   Since $\max E_i <\ga_i \le
\min E_{i+1}$, $i<d-1$, it follows that $\max E_i<\eta_i\le \min E_{i+1}$, $i<d-1$, so $(E_i)_{i<d}$ is $\mc
B$-admissible.
\fprucl
\fprue

Recall that the Cantor-Bendixson index $\varrho_{\mathrm{CB}}(\mc F)$ of a \emph{compact} family $\mc F$ of
finite subsets of $\ka$ is its Cantor-Bendixson index of the family viewed as a closed subset of $2^{\ka}$.
Since $\mc F$ is a scattered compactum, there is an ordinal $\al$ such that the
$\al^{\mathrm{th}}$-derivative of $\mc F$ is empty. Thus, $\varrho_{\mathrm{CB}}(\mc F)$ is the first $\al$
such that the $\al+1$-derivative of $\mc F$ is empty.

\prop\label{iooiewoiwe} Suppose that $\mc F$ is a compact and hereditary family of subsets of some infinite set $C\con \ka$.
\begin{enumerate}
\item[(a)] If there is $\ga\in C$ and some infinite subset $A\con C\setminus \{\ga\}$ such that $\conj{\{\ga\}\cup s}{s\in [A]^{n-1}}\con \mc
F$, then $\varphi_{\mathrm{CB}}(\mc B)\ge n$.
\item[(b)] If
$\varrho_{\mathrm{CB}}(\mc F)\ge n$, then there is $\{\ga_0<\dots<\ga_{n-1}\}$ in $\mc F$ such that
$]\ga_i,\ga_{i+1}[\cap C\neq\buit$ for every $i<n-1$.
\end{enumerate}
\fprop
\prue
(a): We prove that $\{\ga\}\in \partial^{(n)}\mc F$ by induction on $n$.  Suppose that $\conj{\{\ga\}\cup
s}{s\in [A]^{n}}\con \mc F$. Then $\conj{\{\ga\}\cup s}{s\in [A]^{n-1}}\con \partial(\mc F)$, so by inductive
hypothesis, $\{\ga\}\in \partial^{(n-1)}(\partial \mc F)=\partial^{(n)}(\mc F)$, and we are done.

(b): Induction on $n$. Suppose that $\varrho_{\mathrm{CB}}(\mc F)\ge n+1$. Since
$\varrho_{\mathrm{CB}}(\partial \mc F)\ge n$, there is $s=\{\ga_0<\dots<\ga_{n-1}\}\in \partial \mc F$ such
that $]\ga_i,\ga_{i+1}[\cap C\neq \buit$, $i<n-1$. Since $s\in \partial \mc F$ there is a sequence
$(s_k)_{k\in \N}$  in $\mc F\setminus \{s\}$ with limit $s$. Since $\mc F$ is hereditary, we may assume that
$s_{k}=s\cup \{\eta_k\}$ for some $\ga_k\notin s$, $k\in \N$. Suppose first that there are  $k_0,k_1$ such
that $\eta_{k_0}<\eta_{k_1}<\ga_0$. Then $\{\eta_{k_0}\}\cup s$ is the desired set. Suppose now there are
$k_0,k_1$ such that $\ga_{n-1}<\eta_{k_0}<\eta_{k_1}$. Then $\{\eta_{k_1}\}\cup s$ is the desired set.
Otherwise, there is $i<n-1$ and there are $k_0,k_1,k_2$ such that
$\ga_{i}<\eta_{k_0}<\eta_{k_1}<\eta_{k_2}<\ga_{i+1}$, and then $s\cup \{\eta_{k_1}\}$ is the desired
set.\fprue

\lema\label{5roi4tjhiogf}
Suppose that Player I has a winning strategy $\sig$ for the game $\Game_{\mc B,n}$. Let  $C\con \ka$ be a set
closed under $\sig$ and unbounded in itself. Suppose that $\ga_1<\dots<\ga_{n-1}$ are successor elements of
$C$. Then  $\{\min C,\ga_1,\dots,\ga_{n-1}\}\in \varpi_{C}(\mc B)$. Consequently,
\begin{enumerate}
\item[(a)] $\conj{s\in [C]^n}{\min s=\min C \text{ and $s\setminus \{\min s\}\con C^+$}}\con \varpi_{C}(\mc
B)$, and $\varrho_{\mathrm{CB}}(\mc B)\ge n$, if $\mc B$ is compact.
\item[(b)] $[C^+]^{n-1}\con \varpi_{C}(\mc B)$, and
\item[(c)] every block sequence $(s_i)_{i<d}$ of finite subsets of $C$ of length $d\le n$ is $\mc B$-admissible.
\end{enumerate}

\flema
\prue
Suppose that $\{\ga_1<\dots<\ga_{n-1}\}\con C^+$, and for $0<i<n$, let $\de_i\in C$ be such that
$]\de_i,\ga_i]=\{\ga_i\}$.
 We play the following run in the game
$\Game_{\mc B,n}$: Player $I$ plays $\sig(\buit)$.  Since $C$ is closed under $\sig$, it follows that
$\sig(\buit)<\de_1$. Now let Player $II$ play $\de_1$. Next, Player $I$ plays $\de_1<\sig(\de_1)$. It follows
then that $\sig(\de_1)<\ga_1\le\de_2$, so let Player $II$ play $\de_2$, and so on. At the end of the game we
see that $\de_i<\sig(\de_1,\dots,\de_i)<\ga_i$ for every $i=1,...,n-1$, and since $]\de_i,\ga_i]=\{\ga_i\}$,
it follows that $\varpi_{C}(\sig(\de_1,\dots,\de_i))=\ga_i$ for every $i=1,\dots,n$. Since $\sig$ is winning,
$\{\sig(\buit),\sig(\de_1),\dots,\sig(\de_1,\dots,\de_{n-1})\}\in \mc B$. But
$\{\sig(\buit),\sig(\de_1),\dots,\sig(\de_1,\dots,\de_{n-1})\}=\{\min C,\ga_1,\dots,\ga_{n-1}\}$, and we are
done.

(a) Follows from the previous and Proposition \ref{iooiewoiwe} (a). (b) is a trivial consequence of (a). Let
us prove (c): Let $(s_i)_{i<n}$ be a block sequence of finite subsets of $C$.  By adding convenient ordinal
numbers to the sets $s_i$, $i<n$  we may assume wlog that   $\min s_0=\min C$, and that $\min s_i$ are not
limit elements of $C$.    Hence $\{\min s_i\}_{i<n}\in\varpi_{C}(\mc B)$, and so, by Proposition
\ref{mnkdwewrdfnknke}, $(s_i)_{i<n}$ is $\mc B$-admissible.
\fprue

\lema\label{iweiiofgff}
Suppose that Player II has a winning strategy $\sig$ for the game $\Game_{\mc B,n}$. Let $C\con \ka$ be a set
closed under $\sig$ which is unbounded in itself. Suppose that $\ga_0<\dots<\ga_{n-1}$ are in $C$ and are
such that $]\ga_i,\ga_{i+1}[\cap C\neq \buit$ for every $i<n-1$. Then $\{\ga_i\}_{i<n}\notin \varpi_{C}(\mc
B)$. Consequently,
\begin{enumerate}
\item[(a)] If $D\con C$ is such that $]\ga,\eta[\cap C\neq \buit$ for every $\ga<\eta$ in $D$, then $[D]^{n}\cap \varpi_{C}(\mc
B)=\buit$.
\item[(b)]   $\varpi_{C}(\mc B)\con [C]^{<2n-1}$, and if
$(s_i)_{i<d}$ is a $\mc B$-admissible sequence of finite subsets of $C$ then
\begin{equation}\label{j4ithjri4jri4}
\text{$d<2n-1$ and $|\conj{i<d}{|s_i|\ge 2}|\le n$}.
\end{equation}

\item[(c)]   $\varpi_{C}(\mc B)$ is compact and $\varrho_{\mathrm{CB}}(\mc B)\le n$.
\end{enumerate}
\flema
\prue
Let $\ga_0<\dots<\ga_{n-1}$ be as in the hypothesis of the Lemma, and suppose otherwise that
$\{\ga_i\}_{i<n}\in \varpi_{C}(\mc B)$. Fix then $\{\eta_0<\dots<\eta_{n-1}\}\in \mc B$ such that
$\ga_i=\varpi_{C}(\eta_i)$, $i<n$ and let $\de_i\in C$ be such that $\de_i\in ]\ga_i,\ga_{i+1}[ $,  $i<n-1$.
We play the following run in the game $\Game_{\mc B, n+1}$: Player $I$ plays $\eta_0$. Then Player $II$
replies with $\eta_0<\sig(\eta_0)$. Since $C$ is closed under $\sig$, and since $\eta_0\le \ga_0$, it follows
that $\sig(\eta_0)<\de_0<\ga_1$. Since $\varpi_{C}(\eta_1)=\ga_1$, it follows that $\de_0<\eta_1$, and hence
Player I can play $\eta_1$, and Player II replies $\sig(\eta_0,\eta_1)$. By a similar argument,
$\sig(\eta_0,\eta_1)<\eta_2$, so Player I can play $\eta_2$, and so on. In this way Player I was able to
produce $\{\eta_i\}_{i<n}\in \mc B$ against Player II following its winning strategy $\sig$, a contradiction.

(a) is an easy consequence of the previous.  Let us prove (b): Suppose otherwise that
$\{\ga_0<\dots<\ga_{2n-2}\}\in \varpi_{C}(\mc B)$. Then $\{\ga_{2i}\}_{i<n}\in  \varpi_{C}(\mc B)$ and
$\ga_{2i+1}\in ]\ga_{2i},\ga_{2(i+1)}[$, contradicting the statement in the Lemma. Suppose now that
$(s_i)_{i<d}$  is a $\mc B$-admissible sequence of subsets of $C$, and suppose that $d\ge 2n-1$. Wlog we
assume that $d=2n-1$.  Let $\{\eta_i\}_{i<2n-1}\in \mc B$ be such that $\eta_0\le s_0< \eta_1\le
s_1<\dots<\eta_{2n-2}\le s_{2n-2}$.  Then $\{\varpi_{C}(\eta_i)\}_{i<2n-2}\in \varpi_{C}(\mc B)\cap
[C]^{2n-1}$, contradicting the first part of (b).  Similarly one proves that $\conj{i<d}{|s_i|\ge 2}|\le n$.

(c):  Since $\varpi_{C}(\mc B)\con [C]^{<2n-1}$, and $\varpi_{C}(\mc B)$ is hereditary, the family
$\varpi_C(\mc B)$ is compact. Now $\varrho_{\mathrm{CB}}(\varpi_C(\mc B))\le n$ follows from the main
statement in the Lemma, and Proposition \ref{iooiewoiwe}.
\fprue

\lema \label{njkdjfds}
Suppose that $C$ is $\mc B$-strategically closed. Then
\begin{enumerate}
\item[(a)] $\varpi_C(\mc B)\rest C_0=[C_0]^{<\om}$ if $\al(\mc B)=\om$, and where $C_0=C^+\cup \{\min C\}$.
\item[(b)] $\varrho_{\mathrm{CB}}(\varpi_{C}(\mc B))=\al(\mc B)$, when $\al(\mc B)<\om$.
\end{enumerate}
\flema
\prue
(a) readily follows from Lemma \ref{5roi4tjhiogf}. Let us prove (b): From Lemma \ref{iweiiofgff} (c) we know
that $\varpi_{C}(\mc B)\con [C]^{<2n-1}$, so $\varpi_{C}(\mc B)$ is compact. From Lemma \ref{5roi4tjhiogf}
(a) we obtain  that $\varrho_{\mathrm{CB}}(\varpi_{C}(\mc B))\ge \al(\mc B)$.   And
$\varrho_{\mathrm{CB}}(\varpi_{C}(\mc B))\le \al(\mc B)$ is a consequence of Lemma \ref{iweiiofgff} (c).
\fprue
Now Theorem \ref{io43oi4rgrg} is follows readily from the more informative result.

\teor\label{io43oi4rgr22g} Suppose that $\ka$ is an uncountable regular cardinal, and suppose that $C$ is a
club in $\ka$ which is $\mc B$-strategically closed. Then
\begin{enumerate}
\item[(1)] the sequence $(u_\ga)_{\ga\in C}$ in $T_{\theta,\mc B}$ is $\theta^{-1}$-equivalent to the unit basis of
$\ell_1(\ka)$, if $\al(\mc B)=\om$, or
\item[(2)] $\al(\mc B)<\om$ and  the closed linear span of ${\langle u_\ga \rangle_{\ga\in C}}$ in $T_{\theta,\mc B}$
is $c_0$-saturated if $\theta\cdot\al(\mc B)\le 1$ or $\ell_p$-saturated if $\theta\cdot\al(\mc B)> 1$
and where $p=\log(\al(\mc B))/(\log(\al(\mc B))+\log(\theta))$.
\end{enumerate}
\fteor

\prue
Let $C$ be a club of $\ka$ which is $\mc B$-strategically closed.  Suppose first  that $\al(\mc B)=\om$. Then
from Lemma \ref{njkdjfds} (a) we obtain that every finite block sequence $(s_i)_{i<d}$ of finite subsets of
$C$ is $\mc B$-admissible, so it readily follows that $\nrm{\sum_{\ga\in s}a_\ga u_\ga}_{\theta,\mc B}\ge
\theta\sum_{\ga\in s}|a_\ga|$ for every $s\con C$.

Suppose now that $\al(\mc B)<\om$.  Set $X_C:=\overline{\langle u_\ga \rangle_{\ga\in C}}$. Let $X$ be an
infinite dimensional closed subspace of $X_C$. Since $(u_\ga)_{\ga<\ka}$ is a Schauder basis of
$T_{\theta,\mc B}$, we may assume that $X$ is the closed linear span of a normalized block sequence
$(x_n)_{n\in \N}$ with $x_n\in \langle u_\ga \rangle_{\ga\in C}$, for every $n\in \N$.  Let $D:=\bigcup_{n\in
\N}\supp x_n$.  Then $D$ is unbounded in $D$ and it is $\mc B$-strategically closed. It follows that
$\varrho_{\mathrm{CB}}(\varpi_{D}(\mc B))=\al(\mc B)$. It is proved in \cite{Be-De} (see also \cite{Lo-Ma})
that the separable Tsirelson-like space $T_{\theta,\varpi_{D}(\mc B)}$ is saturated by copies of $c_{0}$ if
$\theta\cdot \al(\mc B)\le 1$, or by copies of $\ell_p$ if $\theta\cdot\al(\mc B)>1$, with $p=\log(\al(\mc
B))/(\log(\theta)+\log(\al(\mc B)))$. But by Proposition \ref{kloihreiumnnds}, $X_C$ and
$T_{\theta,\varpi_{D}(\mc B)}$ are the same space, so we are done.
\fprue
%
%
%

\cor\label{oireihjiejf}
Suppose that $\mc B$ is a large hereditary family on an uncountable regular cardinal $\ka$. Then $\al(\mc
B)=\om$ and consequently, the space $T_{\theta,\mc B}$ contains a copy of $\ell_1(\ka)$.
\fcor
\prue
Let $C$ be unbounded in itself and  $\mc B$-strategically closed. Since  $[C]^n\cap \mc B \neq \buit$, for
every $n$, we obtain $[C]^n\cap \varpi_C(\mc B)\neq \buit$. So, by Lemma \ref{njkdjfds}, $\al(\mc B)=\om$ and
the desired result follows from Theorem \ref{io43oi4rgr22g}.
\fprue

A well-known procedure in producing compact and hereditary families in $\om$ is given by the following
operation: Given two families $\mc B$ and $\mc C$ in $\om$, we define
$$\mc B\otimes \mc C:=\conj{\bigcup_{i<d}s_i}{(s_i)_{i<d}\text{ is a block sequence of elements of $\mc B$ such that $\{\min s_i\}_{i<d}\in  \mc C$} }.$$
It is not difficult to see that $\mc B\otimes \mc C$ is compact if both $\mc B$ and $\mc C$ are compact, and
indeed $\varrho_{\mathrm{CB}}(\mc B\otimes \mc C)=\varrho_{\mathrm{CB}}(\mc B)\cdot\varrho_{\mathrm{CB}}(\mc
C)$. In particular, $[\om]^{\le 2}\otimes\mc S$ is a compact family with Cantor-Bendixson rank $\om$, so
comparable to $\mc S$. Similarly, we can produce the $n$-powers of $S$, $\mc S_n:=\mc S \otimes
\overset{(n)}{\cdots}\otimes\mc S$ with $\varrho_{\mathrm{CB}}(\mc S_n)=\om^n$. In the uncountable case, the
situation is again very much different.

\prop\label{i4jtirffdgdf}
Suppose that $\mc B$ is a large family on an uncountable regular cardinal $\ka$. Then there is some subset
$\Ga$ of $\ka$ of cardinality $\ka$ such that $[\Ga]^\om \con \overline{[\ka]^{\le 2}\otimes\mc B}$.
\fprop
\prue
Let $C$ be a club of $\ka$ which is $\mc B$-strategically closed, and let $A$ be an arbitrary subset of $C$
of cardinality $\ka$ such that $\min A\cap C$ and $]\ga,\eta[\cap C$ are infinite sets for every $\ga<\eta$
in $A$. We claim that $[A]^\om\con \overline{\mc B^2}$:  Let $B\in [A]^\om$,  $\ga_i:=\theta_{\om,B}(i)$ for
every $i<\om$, and for each $i<\om$, let $\{\de_j^i\}_{j<\om}\con C$ be such that
$\de_0^0<\de_1^0<\dots<\de_j^0<\ga_0$, and $\ga_i<\de_0^{i+1}<\de_1^{i+1}<\dots<\de_j^{i+1}<\dots<\ga_{i+1}$
for every $i<\om$. Fix $n<\om$. Since we have seen in Corollary \ref{oireihjiejf} that $\al(\mc B)=\om$,
Lemma \ref{njkdjfds}, implies that
\begin{equation}
\{\de_n^0,\ga_0,\de_n^1,\ga_1,\dots,\de_{n}^n,\ga_n\}\in  \varpi_C(\mc B).
\end{equation}
So there is   $\{\eta_i\}_{i<\le n}\in \mc B$ such that
\begin{equation}
\de_n^0<\eta_0\le \ga_0<\de_{n}^1<\eta_1<\ga_1<\dots<\de_n^i<\eta_i<\ga_i<\dots<\de_n^n<\eta_n<\ga_n.
\end{equation}
Then $s_n:=\{\eta_0,\ga_0,\eta_1,\ga_1,\dots,\eta_n,\ga_n\}\in [\ka]^{\le 2}\otimes\mc B$, and $s_n\to_n
\{\ga_0,\dots,\ga_n\}$. Since $n$ is arbitrary, this implies that $B$ is in the topological closure of
$[\ka]^{\le 2}\otimes \mc B^2$.
\fprue

We shall need the following classical result (see, for example, \cite{Wi}).
\begin{thm*}\footnote{Using the Erd\~os-Rado arrow notation, this result is stated as $\ka\rightarrow(\ka,\om)^2$.}[Erd\"{o}s-Dushnik-Miller]
 For every colouring  $c:[\ka]^2\to 2$ either there is $A\in [\ka]^\ka$ such that
$c\rest [A]^2$ is constant with value 0, or there is $B\in [\ka]^\om$ such that $c\rest [B]^2$ is constant
with value 1.
\end{thm*}

\prop
Suppose that $\mc B$ is a large family on an uncountable cardinal number $\ka$. Then $\mc B^2$ is not
compact. Indeed,
\begin{enumerate}
\item If $\ka$ is regular, there is a subset $A$ of $\ka$ of cardinality $\ka$ such that $[A]^\om \con \overline{\mc
B^2}$.
\item If $\ka$ is singular,  for every cardinal number $\la<\ka$ there is a subset $A_\la$ of $\ka$ of cardinality $\la$ such that $[A_\la]^\om \con
\overline{\mc B^2}$.
\end{enumerate}
\fprop
\prue
The singular case is a direct consequence of the regular case. So, we assume that $\ka$ is an uncountable
regular cardinal number.  Let $c:[\ka]^2\to 2$ defined for $s\in [\ka]^2$ by $c(s)=0$ iff $s\in \mc B$. By
the Erd\"{o}s-Dushnik-Miller theorem there must be a subset $A\con \ka$ of cardinality $\ka$ such that $[A]^2\con
\mc B$, since the other possibility would imply that $\mc B$ is not large. Hence $[A]^{\le }\otimes \mc
B\rest A\con (\mc B\rest A)^2$, and the desired result follows from Proposition \ref{i4jtirffdgdf}.
\fprue

\nota
 Recall that a family $\mc B$ is called \emph{spreading} when if $(\{\ga_i\})_{i<d}$ is $\mc
B$-admissible, then $\{\ga_i\}_{i<d}\in \mc B$.  In particular, $\varpi_C(\mc B)=\mc B\rest C$ for every
$C\con \ka$ unbounded in $C$.  So, if $\mc B$ is a   hereditary and spreading family on an uncountable
cardinal $\ka$, then there is an uncountable subset $C$ of $\ka$ and some ordinal $\al\le \om$ such that $\mc
B\rest C=[C]^{<\al}$:  This is basically consequence of Lemma \ref{5roi4tjhiogf} and Lemma \ref{njkdjfds}.
\fnota

\nota
It is a well-known fact that the separable Tsirelson space $T$ does not contain isomorphic copies of $c_0$ or
$\ell_p$, $p\ge 1$. Indeed, $T$ does not have subsymmetric basic sequences.  One may ask if there is some
cardinal number $\ka$ such that every Banach space  of such density $\ka$ contains an isomorphic copy of
either $c_0$ or $\ell_p$, $p\ge 1$.  This is not true, as the following simple remark shows: Let $X$ be the
separable space constructed by Figiel and Johnson \cite{Fi-Jo} with a subsymmetric\footnote{indeed,
symmetric} basis $(v_n)_{n\in \om}$ and not containing $c_0$ or $\ell_p$, $p\ge 1$. Then define on
$c_{00}(\ka)$ the following norm: For $x\in c_{00}(\ka)$, let
\begin{equation}\label{jknbkjgnjg}
\nrm{x}_{X,\ka}:=\nrm{\theta_{\supp x, |\supp x|}(x)}_X
\end{equation}
where $\theta_{\supp x, |\supp x|}(x)=\sum_{\ga\in \supp x} (x)_\ga v_{\theta_{\supp x,|\supp x|}(\ga)}.$ The
formula in \eqref{jknbkjgnjg} defines a norm on $c_{00}(\ka)$ because $(v_n)_{n}$ is subsymmetric. Moreover,
$(u_\ga)_{\ga<\ka}$ is a subsymmetric basis of the completion $\mathfrak{X}_{X,\ka}$ of
$(c_{00}(\ka),\nrm{\cdot}_{X,\ka})$. A standard gliding-hump argument shows that $\mathfrak{X}_{X,\ka}$ does
not contains $c_0$ or $\ell_p$, $p\ge 1$.
\fnota


\end{document}